\documentclass[10pt]{amsart}

\usepackage{amsmath,amssymb,amsfonts,epsfig}

\newcommand{\rb}{\raisebox}
\newcommand{\ig}{\includegraphics}
\newcommand\risS[6]{\rb{#1pt}[#5pt][#6pt]{\begin{picture}(#4,15)(0,0)
  \put(0,0){\ig[width=#4pt]{#2.eps}} #3
     \end{picture}}}
\def\a{\alpha}
\def\b{\beta}
\def\d{\delta}
\def\cS{\mathcal S}
\def\bb{\mathbf b}
\def\bK{\mathbf K}
\newtheorem{thm}{Theorem}[section]
\newtheorem{defn}[thm]{Definition}
\newtheorem{exa}[thm]{Example}

\newtheorem{lemma}[thm]{Lemma}
\newtheorem{prop}[thm]{Proposition}

\newtheorem{rem}[thm]{Remark}

\def\kb#1{[ #1 ]}

\def\wo{\overline}

\def\wt{\widetilde}

\def\pn#1#2#3{\put(#1,#2){\mbox{\tt\scriptsize #3}}}

\def\smiinc{\mbox{\begin{picture}(5,5)(0,0)
          \put(3,2){\circle{7}} \put(-.2,.1){\mbox{\scriptsize $-$}}
                 \end{picture}}}
\begin{document}

\title[Arrow ribbon graphs]{Arrow ribbon graphs}
\author{ROBERT BRADFORD, CLARK BUTLER, SERGEI CHMUTOV}
\date{}

\address{Department of Mathematics, University of Kansas,
1460 Jayhawk Blvd., Lawrence, Kansas 66045\linebreak 
{\tt rbradford@math.ku.edu}\linebreak 
Department of Mathematics, The Ohio State University,
231 West 18th Avenue, Columbus, OH 43210
{\tt butler.552@buckeyemail.osu.edu}\qquad
{\tt chmutov@math.ohio-state.edu}}


\begin{abstract}
We introduce an additional arrow structure on ribbon graphs.
We extend the dichromatic polynomial to ribbon graphs with this structure. This extended polynomial satisfies the contraction-deletion relations and behaves naturally with respect to the partial duality of ribbon graphs. 
From a virtual link, we construct an arrow ribbon graph whose extended dichromatic polynomial specializes to the arrow polynomial of the virtual link recently introduced by H.~Dye and L.~Kauffman. This result generalizes the classical Thistlethwaite theorem to the arrow polynomial
of virtual links.
\end{abstract}

\maketitle

\section*{Introduction} \label{s:intro}
The classical Thistlethwaite theorem \cite{Th} relates the Jones polynomial $V_L(t)$ of an alternating link $L$ to the Tutte polynomial 
$T_{G_L}$ of an appropriate planar graph $G_L$
$$\risS{-7}{kd-gr}{\put(3,40){\mbox{$L$}}  \put(187,40){\mbox{$G_L$}}
\put(-110,25){\mbox{$\begin{array}{l}
               V_L(t)\ =\ t+t^3-t^4 \vspace{8pt}\\
               \ =\ -t^2(-t^{-1}-t+t^2)
                     \end{array}$}}
\put(205,30){\mbox{$\begin{array}{l}
               T_{G_L}(x,y)=y+x+x^2 \vspace{8pt}\\
               T_{G_L}(-t,-t^{-1})= \\
                 \hspace{10pt}=\ -t^{-1}-t+t^2
                     \end{array}$}}
                 }{200}{45}{10} 
$$
For non-alternating links Thistlethwaite \cite{Th} imposed an additional structure on planar graphs $G_L$.
For every crossing that contradicted the alternating pattern, he assigned a sign ``$-$" to the corresponding edge of $G_L$. Thus he worked with signed graphs. 
In 1989 L.~Kauffman  \cite{Ka2} reformulated and reproved the Thistlethwaite theorem by extending the Tutte polynomial to signed graphs. He also extended Thistlethwaite's contraction-deletion property of the Jones polynomial to the signed Tutte polynomial \cite{Th}.
 
The Thistlethwaite theorem was later generalized to virtual links using ribbon graphs (see  \cite{Ch, CP, ChVo}), and also using the relative Tutte polynomial of planar graphs  \cite{DH}. A relation between these two approaches was later established by two of the authors \cite{BuCh}.

Recently, it was observed  \cite{DK,Mi} that for virtual links the Jones polynomial can be split into several parts which are invariant under the Reidemeister moves individually. These additional parts do not arise in the case of classical links. The generating function of these parts was called the {\it arrow polynomial} in  \cite{DK}.

We introduce an additional arrow structure on ribbon graphs which is
inspired by the arrow structure on virtual links used to formulate the
arrow polynomial. This structure enables us to formulate an extension
of Thistlethwaite's theorem analagous to Kauffman's extension \cite{Ka2}. We
obtain the arrow polynomial as a specialization of an appropriate
extension of the Bollobas--Riordan polynomial to ribbon graphs with an
arrow structure. This extended polynomial satisfies
contraction-deletion relations and behaves well with respect to the
partial duality of ribbon graphs \cite{Ch}.

An arrow structure is a choice of arrows tangent to the boundaries of the vertex-discs and the edge-ribbons. Particular cases of arrow structure have appeared in the literature before.
Under the name {\it arrow presentation} it has been used to encode ribbon graphs \cite{Ch,EMM,Mo1,Mo2,Mo3}. A similar structure appeared in the theory of Vassiliev knot invariants  \cite{BN} under the name {\it marked surfaces} which came from Penner's triangulation of the decorated moduli space of Riemann surfaces  \cite{Pe}, see also \cite[Sec.4.4]{LZ}.

This work was done as part of the Summer 2010 undergraduate research working group
\begin{center}\verb#http://www.math.ohio-state.edu/~chmutov/wor-gr-su10/wor-gr.htm#
\end{center}
``Knots and Graphs" 
at the Ohio State University. We are grateful to all participants of the group for valuable discussions, to the OSU Honors Program Research Fund for the student financial support, to Ilya Kofman and Iain Moffatt for useful comments, and to an anonymous referee for various suggestions which drastically improved the exposition of the paper.

\section{Ribbon graphs and arrow structure} \label{s:rg-ast}

By a {\it ribbon graph} we mean an abstract (not necessarily orientable) surface with boundary decomposed into topological discs of two types, 
{\it vertex-discs} and {\it edge-ribbons}, satisfying the following natural conditions:  the vertex-discs and the edge-ribbons intersect by disjoint line segments, each such line segment lies on the boundary of precisely one vertex and precisely one edge, and every edge contains exactly two such line segments. 
We refer to  \cite{BR3,Ch} for precise definitions and to
 \cite{GT,LZ,MT} for the general notions and terminology of topological graph theory. Ribbon graphs are considered up to homeomorphisms of the underlying surfaces preserving the decomposition. A ribbon graph can be regarded as a regular neighborhood of a graph cellularly embedded into a surface. Thus the language of ribbon graphs is essentially the same as for cellularly embedded graphs. Here are a few examples of ribbon graphs.
 $$
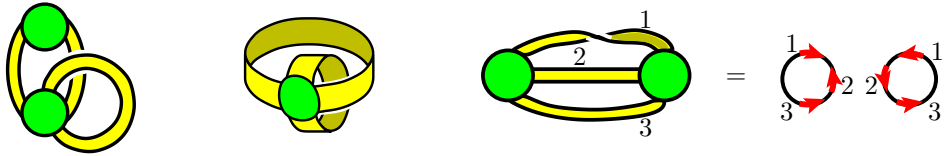
\risS{-27}{rg-ex1}{}{50}{0}{0}\hspace{1.4cm}
  \risS{-19}{rg-ex2}{}{50}{0}{0}\hspace{1.4cm}
\risS{-15}{du-ex33-na}{\put(60,36){$1$}\put(35,22){$2$}
                           \put(60,-5){$3$}}{80}{30}{30}\quad=\quad
\risS{-12}{du-ex33a}{\put(0,-2){\mbox{$3$}}\put(2,24){\mbox{$1$}} 
         \put(23,8){\mbox{$2$}}\put(56,-2){\mbox{$3$}}
         \put(57,21){\mbox{$1$}}\put(32,8){\mbox{$2$}}}{60}{0}{0}
$$
Alternatively, a ribbon graph may be given by an \emph{arrow presentation} (see the third picture). An arrow presentation consists of a set of disjoint circles together with a collection of arrow markings on these circles. These arrows are labeled in pairs. To obtain a ribbon graph from an arrow presentation, we glue discs to each of the circles and attach edge ribbons to each pair of arrows according to the orientation of the arrows.  

\subsection{Arrow structure}
\begin{defn}\rm
An {\it arrow ribbon graph} is a ribbon graph together with
a set (possibly empty) of arrows tangent to the boundaries of the (vertex- and edge-) discs of the decomposition.
Two arrow graphs are {\it equivalent} if there is a homeomorphism between the corresponding surfaces respecting the decompositions and the orientations of the arrows.
\end{defn}
The endpoints of segments along which the edges are attached to the vertices
divide the boundaries of (vertex- and edge-) discs into arcs.

We will refer to the sides of each edge connecting to the vertices as the {\it attaching arcs}. We will refer to the sides which do not attach to vertices as the {\it free edge arcs}. Similarly, we will refer to those arcs on vertices which are the complement of the attaching arcs as {\it free vertex arcs}. 

The arrows may be slid along these arcs by an appropriate homeomorphism, but may not be slid over the end-points of the segments, and so may not change the type of arc that it is on. Each arc may contain several arrows.

For example, if there are several arrows on an attaching arc, then only their order on the arc is relevant, not their actual position on the arc. However, they all must be located on the arc and should not be slid to a free vertex arc or a free edge arc. 

An important example of an arrow structure is the arrow structure given by a particular arrow presentation of a ribbon graph. In this case the arrows all lie on the attaching arcs of the edges. Each pair of arrows corresponds to a map attaching the corresponding edge ribbon to the vertex discs. 


\subsection{Partial duality}

Partial duality of ribbon graphs was introduced in  \cite{Ch} under the name
{\it generalized duality}. Dan Archdeacon suggested a more appropriate term, {\it partial duality}. Under this name it was then used in papers  \cite{EMM,Mo1,Mo2,Mo3,VT}. 

For any ribbon graph $G$, there is a {\it natural} dual ribbon graph $G^*$, also called the
{\it Euler-Poincar\'e} dual. First we glue a disc, also known as a {\it face}, to each boundary component of $G$, obtaining a closed surface $\wt{G}$ without boundary. 
Then we remove the interior of all vertex-discs of $G$.
The newly glued disc-faces will be the vertex-discs of $G^*$. The edge-ribbons for $G^*$ will be the same as for $G$ but now they are attached to the new vertices by the pair of opposite arcs which used to be free edge arcs in $G$, and the attaching arcs of the ribbons become free edge arcs. This is a particular case of the partial duality with respect to the set of all edges of $G$.

\begin{defn}\rm
{\it Partial duality} is duality with respect to a subset 
$D\subseteq E(G)$ of edges of $G$. We denote this partially dual graph by $G^D$. It is constructed as follows. Consider the spanning subgraph $F_D$ of $G$ containing all vertices of $G$ and only the edges from the  subset $D$. 
Regard the boundary components of $F_D$ as curves on the surface of $G$ via the inclusion $F_D \hookrightarrow G$.
Glue a disc to $G$ along each connected component of this curve and remove the interior of all vertices of $G$. Regard these newly glued discs as vertices. 

The result $G^D$ is easily seen to be a ribbon graph: Each edge in $D$ is now attached to a vertex by the pair of opposite arcs which were free edge arcs in $G$, while the attaching arcs are now free edge arcs, since the interiors of the vertices were deleted. The edges of $E \backslash D$ are attached by the same pair of opposite sides as before.  

For arrow ribbon graphs we preserve all of the arrows on the arcs of the edges and vertices. Note however that arrows on the free edge arcs of an edge in $D$ will become arrows on the attaching arcs of that edge in $G^D$, and vice versa.
\end{defn}

The next table illustrates the partial duality with respect to a single edge. 
$$\label{tab:pd}
\begin{array}{|c||c||c|} \hline
e&G & G^e \makebox(0,10){} \\ \hline\hline

\rb{10pt}{\small non loop}
&\risS{-14}{cdpd1}{\put(25,37){$\alpha$}\put(25,0){$\beta$}
             \put(46,26){$\varepsilon$}
             \put(10,22){$\gamma$}\put(47,13){$\delta$}}{60}{45}{25}&
\risS{-14}{cdpd4}{\put(35,29){$\alpha$}\put(35,18){$\beta$}
             \put(25,28){$\varepsilon$}
             \put(2,22){$\gamma$}\put(25,16){$\delta$}}{60}{0}{0}  \\ \hline

\rb{-5pt}{\small orientable loop}
&\risS{-35}{orl0}{\put(40,40){$\alpha$}\put(33,30){$\beta$}
             \put(13,12){$\gamma$}\put(13,52){$\delta$}}{70}{40}{40}&
\risS{-12}{orl1}{\put(32,24){$\alpha$}\put(47,8){$\beta$}
             \put(18,6){$\gamma$}\put(19,13){$\delta$}}{70}{0}{0}\\ \hline

\rb{-5pt}{\small non-orientable loop}
&\risS{-35}{norl0}{\put(40,40){$\alpha$}\put(40,25){$\beta$}
             \put(13,12){$\gamma$}\put(13,52){$\delta$}}{70}{40}{40}&
\risS{-35}{norl1}{\put(23,40){$\gamma$}\put(23,25){$\delta$}
             \put(50,12){$\beta$}\put(50,52){$\alpha$}}{70}{0}{0}\\ \hline
\end{array} 
$$

We label the arrows by letters $\alpha, \beta, \gamma, \delta, \varepsilon$ in order to make clear which arrow goes to which under the partial duality. The boxes \fbox{$A$} and \fbox{$B$} here stand in for the presence of other edges which may be attached to the vertex along the dotted arcs. The order of attachment of these edges in the box \rb{6pt}{\fbox{\rotatebox{180}{$A$}}} is opposite to the one in \fbox{$A$}. 

Here are some more examples (the details of partial duality are worked out in  \cite{Ch}).
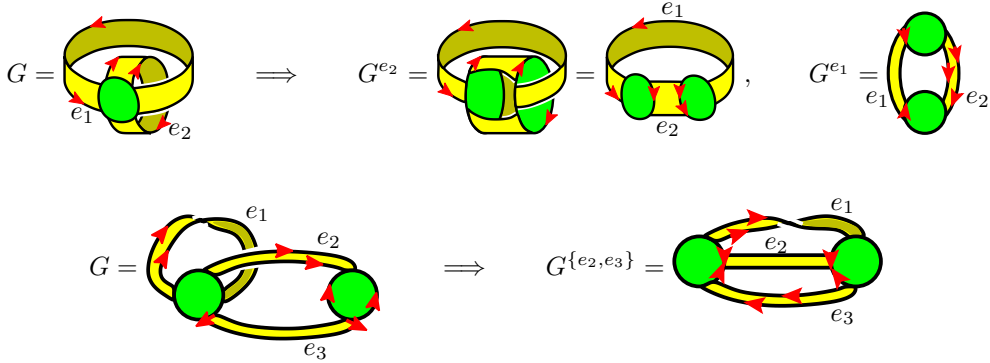
$$G = \risS{-20}{ag-ex2}{\put(3,6){$e_1$}\put(40,0){$e_2$}}{50}{25}{30} \qquad\Longrightarrow\qquad
G^{e_2} = \risS{-20}{du-ex2c0}{}{50}{0}{0} = 
  \risS{-15}{du-ex2c}{\put(22,41){$e_1$}\put(20,-4){$e_2$}}{50}{0}{0}
\ , \qquad 
G^{e_1} = \risS{-20}{du-ex2b}{\put(-8,12){$e_1$}\put(30,12){$e_2$}}{30}{0}{0}
$$
$$G = \risS{-25}{ag-ex3}{\put(38,47){$e_1$}\put(65,37){$e_2$}
                           \put(60,-5){$e_3$}}{90}{30}{30} \qquad\Longrightarrow\qquad
G^{\{e_2,e_3\}} = \risS{-15}{du-ex33c}{\put(60,40){$e_1$}\put(35,25){$e_2$}
                           \put(60,-2){$e_3$}}{80}{0}{0}
$$

{\bf Properties}  \cite{Ch}.
\begin{itemize}
\item[(a)] $G^\emptyset = G$.
\item[(b)] $G^{E(G)} = G^*$.
\item[(c)] $\bigl(G^D\bigr)^{D'}= G^{(D\cup D')\setminus (D\cap D')}$,
 in particular
 \begin{itemize}
   \item $\bigl(G^D\bigr)^D= G$,
   \item for $e\not\in D$, 
$G^{D\cup\{e\}} = \bigl(G^D\bigr)^{\{e\}} = \bigl(G^{\{e\}}\bigr)^D$.
 \end{itemize}
\item[(d)] Partial duality preserves orientability.
\item[(e)] Partial  duality preserves the number of connected components.
\end{itemize}

\subsection{Contraction-Deletion}
\begin{defn}\label{def:contraction}\rm
For an arrow ribbon graph $G$ with an edge $e$, {\it deletion} 
of the edge $e$ gives the arrow ribbon graph $G-e$ obtained from $G$ by removing the edge-ribbon $e$. If the arrows on $e$ were chosen on the attaching arcs of $e$, then we keep those arrows and consider them as arrows on the corresponding free vertex arcs of $G-e$. If the arrows of $e$ were chosen on the free edge arcs of $e$ forming two arcs of the boundary of $G$, then we remove these along with the edge $e$.
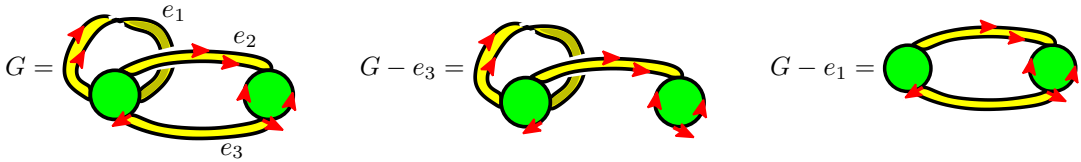
$$G = \risS{-25}{ag-ex3}{\put(38,47){$e_1$}\put(65,37){$e_2$}
                           \put(60,-5){$e_3$}}{90}{30}{35}\hspace{.8cm}
G-e_3 = \risS{-25}{c-e3}{}{90}{0}{0}\hspace{.8cm}
G-e_1 = \risS{-14}{c-e1}{}{75}{0}{0}
$$
\end{defn}

\begin{defn}\rm
{\it Contraction} is defined using partial duality. We define the contraction of an edge $e$ in a ribbon graph $G$ by 
$$G/e\ :=\ G^{\{e\}}-e\ .
$$
\end{defn}

This definition coincides with the usual notion of contraction for ribbon graphs (For details, see  \cite{Ch}). We have the following properties: for $e\not\in D$,

$$(G/e)^D = G^D/e = G^{D\cup e}-e\qquad\mbox{and}\qquad
(G-e)^D = G^D-e = G^{D\cup e}/e\ .
$$

Since we have specified how the arrow structure behaves under deletion of edges and partial duality, we have also specified how our arrow structure behaves under contraction of edges. We present here the particular cases of contraction and deletion of a non-loop, an orientable loop, and a non-orientable loop.




$$\label{tab:cont-del-pd}
\begin{array}{|c||c||c||c|c|} \hline
e&G & G^e & G-e=G^e/e & G/e=G^e-e\makebox(0,10){} \\ \hline\hline

\rb{10pt}{\small non loop}
&\risS{-14}{cdpd1}{\put(25,37){$\alpha$}\put(25,0){$\beta$}
             \put(46,26){$\varepsilon$}
             \put(10,22){$\gamma$}\put(47,13){$\delta$}}{60}{45}{25}&
\risS{-14}{cdpd4}{\put(35,29){$\alpha$}\put(35,18){$\beta$}
             \put(25,28){$\varepsilon$}
             \put(2,22){$\gamma$}\put(25,16){$\delta$}}{60}{0}{0}&
\risS{-14}{cdpd2}{\put(46,26){$\varepsilon$}
             \put(10,22){$\gamma$}\put(47,13){$\delta$}}{60}{0}{0}
&\risS{-14}{cdpd3}{\put(25,27){$\alpha$}\put(25,15){$\beta$}}{60}{0}{0}  \\ \hline

\rb{-5pt}{\small orientable loop}
&\risS{-35}{orl0}{\put(40,40){$\alpha$}\put(33,30){$\beta$}
             \put(13,12){$\gamma$}\put(13,52){$\delta$}}{70}{40}{40}
&\risS{-12}{orl1}{\put(32,24){$\alpha$}\put(47,8){$\beta$}
             \put(18,6){$\gamma$}\put(19,13){$\delta$}}{70}{0}{0}&
\risS{-15}{orl2}{\put(40,20){$\alpha$}}{70}{0}{0}&
\risS{-12}{orl3}{\put(47,8){$\beta$}
             \put(18,6){$\gamma$}\put(19,13){$\delta$}}{70}{0}{0}\\ \hline

\rb{-5pt}{\small non-orientable loop}
&\risS{-35}{norl0}{\put(40,40){$\alpha$}\put(40,25){$\beta$}
             \put(13,12){$\gamma$}\put(13,52){$\delta$}}{70}{40}{40}&
\risS{-35}{norl1}{\put(23,40){$\gamma$}\put(23,25){$\delta$}
             \put(50,12){$\beta$}\put(50,52){$\alpha$}}{70}{0}{0}&
\risS{-18}{norl2}{\put(40,23){$\alpha$}\put(40,8){$\beta$}}{70}{0}{0}&
\risS{-18}{norl3}{\put(23,21){$\gamma$}\put(23,8){$\delta$}}{70}{0}{0}\\ \hline
\end{array} 
$$

Here are three more examples. 
$$\begin{array}{r@{\qquad\Longrightarrow\qquad}l}
G = \risS{-20}{ag-ex2}{\put(3,6){$e_1$}\put(40,0){$e_2$}}{50}{30}{30} &
G/e_2 = \risS{-15}{b-mod-e}{}{50}{0}{0}\\

G = \risS{-25}{ag-ex3}{\put(38,47){$e_1$}\put(65,37){$e_2$}
                           \put(60,-5){$e_3$}}{90}{30}{35} &
G/e_3 = \risS{-25}{c-mod-e3}{}{50}{0}{0}\hspace{.8cm}
G/e_1 = \risS{-14}{c-mod-e1}{}{75}{0}{0} 
\end{array}
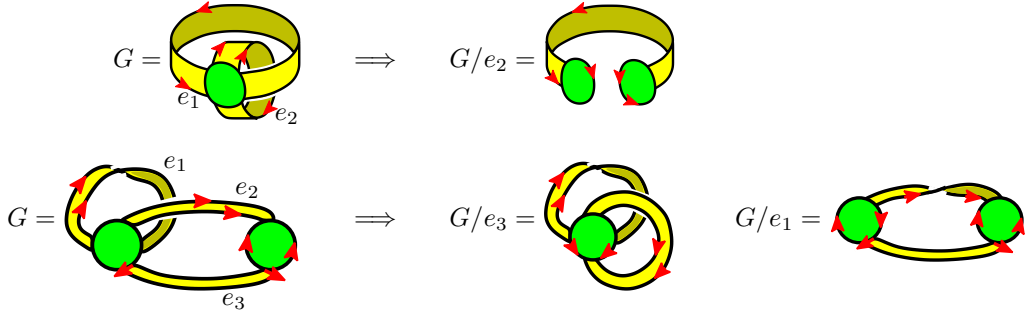
$$

\section{Arrow dichromatic polynomial} \label{s:ar-BR}

Tutte's dichromatic polynomial $Z_G(a,b)$, also known as a partition function of the Potts model in statistical mechanics, was generalized to signed graphs in  \cite{Ka2}. Its multivariable version was introduced in \cite{Tr} and used in \cite{Sok}. The multivariable Tutte polynomial also appears as a very special case of Zaslavsky's colored Tutte polynomial of a matriod \cite{Za} 
and consequently also of Bollob\'as and Riordan's colored Tutte
polynomial of graph \cite{BR1}.
It can be defined as 
$$Z_G(a,\bb):=\sum_{F\subseteq E(G)} a^{k(F)} \prod_{e\in F}b_e\ ,
$$
The sum runs over all spanning subgraphs of $G$, which we identify with subsets $F$ of $E(G)$. $\bb:=\{b_e\}$ is the set of variables (weights) $b_e$ corresponding to the edges $e$ of $G$, and $k(F)$ denotes the number of connected components of $F$.

The multivariable dichromatic polynomial was generalized to ribbon graphs in  \cite{Mo} as
$$Z_G(a,\bb,c):=\sum_{F\subseteq E(G)} a^{k(F)} 
   \Bigl(\prod_{e\in F}b_e\Bigr) c^{bc(F)}\ ,
$$
where $bc(F)$ is the number of connected components of the boundary of $F$. Its signed version from  \cite{VT} can be obtained by substitution
$$\begin{array}{rcl}
a&=&q,\\
b_e&=&\left\{\begin{array}{ll}
             \a_e&\mbox{if $e$ is positive,}\\
             q/\a_e&\mbox{if $e$ is negative,}
           \end{array}\right.
\end{array}
$$
and multiplication of the whole polynomial by 
$\displaystyle \prod_{e\in E(G)} q^{-1/2}\a_e$.

\begin{rem}\label{rem:mBR} \rm
The dichromatic polynomial is essentially equivalent to the Tutte polynomial.
Its ribbon graph formulation, known as the Bollob\'as-Riordan polynomial  \cite{BR3}, is equivalent to the ribbon graph  formulation of the dichromatic polynomial in the same way. 
Similarly, one may introduce a multivariable Bollob\'as-Riordan polynomial \cite{Mo,VT}. Sometimes it is more convenient to use a homogeneous, doubly weighted form of it, which can be defined as
\begin{equation}\label{eq:mBR}
BR_G(X,Y,Z):=\sum_{F\subseteq E(G)} \ 
  (\prod_{e\in F}x_e)\ (\prod_{e\in E(G)\setminus F}y_e)\  
  X^{r(G)-r(F)}\ Y^{n(F)}\ Z^{k(F)-bc(F)+n(F)}\ ,
\end{equation}
where $r(F):=|V(G)|-k(F)$ is the {\it rank} of $F$, $n(F):=|E(F)|-r(F)$ is the {\it nullity} of $F$, and with each edge $e$ we associate a pair of variables $(x_e,y_e)$.

Of course, $BR_G(X,Y,Z)$ is equivalent to $Z_G(a,\bb,c)$ due to the relation
$$BR_G(X,Y,Z)= \Bigl(\prod_{e\in E(G)}y_e\Bigr) (YZ)^{-v(G)}X^{-k(G)}
Z_G(XYZ^2,\{x_eYZ/y_e\},Z^{-1})\ .
$$

A signed version of the Bollob\'as-Riordan polynomial, which was introduced in  \cite{CP} 
and used in \cite {ChVo,Ch}, can be obtained from the multivariable Bollob\'as-Riordan polynomial 
by choosing the weights $(x_+,y_+)$  (resp. $(x_-,y_-)$) of positive (resp. negative) edges to be
$$x_+:=y_+:=1,\qquad\qquad 
x_-:=\sqrt{\frac{X}{Y}},\ y_-:=\sqrt{\frac{Y}{X}}\ .
$$
The main combinatorial results of  \cite{Ch} about contraction-deletion and partial duality may be generalized to the doubly weighted Bollob\'as-Riordan polynomial in a straightforward way. 
\end{rem}

\begin{defn}\label{def:aBR}\rm
In the presence of an arrow structure we can extend the dichromatic polynomial $Z_G(a,\bb,c)$ to the {\it arrow dichromatic polynomial}, 
$$A_G(a,\bb,c,\bK):=\sum_{F\subseteq E(G)} a^{k(F)} 
   \Bigl(\prod_{e\in F}b_e\Bigr) c^{bc(F)}
   \prod_{f\in\partial(F)}K_{i(f)}\ ,
$$
where $F$ is a spanning subgraph of $G$ which we will also refer to as a 
{\it state}; the parameters $k(F)$ and $bc(F)$ are the same as before; and the rightmost product runs over all boundary components $f$ of $F$. The variables $K_{i(f)}$ are assigned to each boundary component $f$ according to the arrangement of arrows along this boundary component. Namely, the subscript $i(f)$ is equal to half of the number of arrows along the boundary component $f$ remaining after recursive cancellations of all neighboring pairs of arrows which point in the same direction:
\end{defn}
$$\risS{-25}{ar-can1}{}{50}{20}{30}\quad\risS{-2}{toto}{}{25}{0}{0}\ 
      K_1\qquad
\risS{-25}{ar-can2}{}{50}{20}{30}\quad\risS{-2}{toto}{}{25}{0}{0}\ 
      K_{1/2}\qquad
\risS{-25}{ar-can3}{}{50}{20}{30}\quad\risS{-2}{toto}{}{25}{0}{0}\ 
      K_2
$$
We set $K_0=1$. Note that whenever the number of arrows is odd on a boundary component, the associated variable is always $K_{1/2}$.
The arrow dichromatic polynomial is a polynomial in infinitely many variables 
$a, b_e, c, K_{1/2}, K_1, K_2, \dots$.
However, for a concrete graph $G$ only finitely many $K$'s appear in $A_G(a,\bb,c,\bK)$.

\begin{exa}\rm
For the arrow graph $G$ shown on the leftmost column in the table, there are eight states. Their parameters and the corresponding monomial in $K$'s are shown.
$$\begin{array}{c||c|c|c|c} \label{br-table}
\!\!\risS{8}{ag-ex3}{\put(28,40){$e_1$}\put(50,32){$e_2$}
                           \put(40,7){$e_3$}}{75}{55}{0}
& \risS{8}{agBBB}{}{70}{0}{0}
& \risS{8}{agBBA}{}{70}{0}{0} 
& \risS{15}{agBAB}{}{70}{0}{0}  
& \risS{15}{agBAA}{}{70}{0}{0}\\ \hline
 k,bc,\prod K_{i(f)} & 1,2,K_1^2 & 1,1,K_1 & 1,1,K_1 & 2,2,K_{1/2}^2
        \makebox(0,12){}\\ \hline\hline
& \!\risS{8}{agABB}{}{75}{55}{0}\!
& \!\risS{8}{agABA}{}{75}{0}{0}\!
& \!\risS{10}{agAAB}{}{75}{0}{0}\!
& \!\risS{10}{agAAA}{}{75}{0}{0}\!\!\!\!\!\\ \cline{2-5}
& 1,1,K_1 & 1,1,K_1 & 1,1,1 & 2,2,K_{1/2}^2 \makebox(0,12){} 
\end{array}
$$
Thus
$$
A_G = ab_2b_3c^2K_1^2+ab_3cK_1+ab_2cK_1+
              a^2c^2K_{1/2}^2 
 +ab_1b_2b_3cK_1+ab_1b_3cK_1+ab_1b_2c+a^2b_1c^2K_{1/2}^2\ .
$$
\end{exa}

\begin{rem}\rm 
We do not have to make the cancellation of arrows in the definition above to obtain our combinatorial results in this section. We can simply treat the 
arrangement of arrows on a circle $f$ as a formal variable of $K$-type. However, we will need this cancellation for the arrow generalization of the Thistlethwaite theorem, 
where it corresponds to invariance of the arrow polynomial of  \cite{DK} under the Reidemeister moves.
\end{rem}

Now we are ready to formulate the contraction-deletion properties of the arrow polynomial. 

\begin{prop} {\bf The contraction-deletion properties.}\\
The arrow dichromatic polynomial $A_G(a,\bb,c,\bK)$ possesses the following properties.
\begin{align*}A_{G_1\sqcup G_2} &= A_{G_1}\cdot A_{G_2}\ ;\\
A_G &= \left\{\begin{array}{ll}
A_{G-e} + b_eA_{G/e} & \mbox{if $e$ is not an orientable loop,}\\
A_{G-e} + (b_e/a)A_{G/e} & \mbox{if $e$ is a trivial orientable loop.}
\end{array}\right.
\end{align*}
\end{prop}

\begin{proof}
The first property of multiplicativity under the disjoint union 
$G_1\sqcup G_2$ is obvious. 
The proof of the contraction-deletion properties follows the standard procedure. 
One can split the set of spanning subgraphs $F$ of $G$ into two types according to the property $e\in F$ or $e\not\in F$. The subgraphs of the first (resp. second) type may be regarded as spanning subgraphs of $G/e$ (resp. $G-e$). 
For an edge $e$ which is not an orientable loop, the exponents of variables $a$ and $c$ will be preserved when we consider $F$ as a subgraph of $G$, or of either $G/e$ or $G-e$. Also the corresponding monomials in $K$'s will be equal as one can see from the table on page \pageref{tab:cont-del-pd}.
This implies the first contraction-deletion property.
The second contraction-deletion property follows from the fact that for
a trivial orientable loop $e$, a state of $G/e$ corresponding to a subgraph $F\ni e$ of $G$ always has one more connected component than $F$, i.e. $k(F)$ increases by 1 when we are passing to the contraction $G/e$ according to our definition \ref{def:contraction} of the contraction of a loop.
\end{proof}

\begin{rem} \rm
There is no general contraction-deletion property for a non-trivial orientable loop. However, as in \cite[Lemma 3.3]{Ch}, we have such a property for an evaluation of the arrow dichromatic polynomial at $a=1$. Namely, for any edge $e$:
$$A_G(1,\bb,c,\bK)=A_{G-e}(1,\bb_{\not=e},c,\bK) + 
b_eA_{G/e}(1,\bb_{\not=e},c,\bK)\ ,
$$
where $\bb_{\not=e}=\{b_{e'}\}_{e'\in E(G)\setminus e}$.
For a non-trivial orientable loop $e$ it follows from the partial duality below.
\end{rem}

\begin{prop} \label{prop:pardu}{\bf The partial duality properties.}\\
Let $D\subseteq E(G)$ be a subset of edges and $G':=G^D$ be the corresponding partial dual arrow graph.
The evaluation of the arrow dichromatic polynomial at $a=1$ satisfies to the equation:
$$A_G(1,\bb,c,\bK) = \Bigl(\prod_{e\in D} b_e\Bigr)
  A_{G'}(1,\bb_D,c,\bK)\ , 
$$
where the weights $\bb_D=\{b'_e\}$ of edges of $G'$ are 
$$b'_e = \left\{\begin{array}{ll}
b_e & \mbox{if $e\not\in D$\ ,}\\
1/b_e & \mbox{if $e\in D$\ .}
\end{array}\right.
$$
\end{prop}

\begin{proof} The proof is similar to \cite[Theorem 3.1]{Ch}.
The 1-to-1 correspondence between the spanning subgraphs $F$ of $G$ and the spanning subgraphs $F'$ of $G'$ is given by the symmetric difference:
$F'=F\Delta D:=(F\cup D)\setminus (F\cap D)$.

This correspondence assures that the monomials in weights $b_e$ are equal to each other for $F$ and $F'$. Indeed,
$$\Bigl(\prod_{e\in D} b_e\Bigr) \prod_{e'\in F'} b'_e =
\Bigl(\prod_{e\in D} b_e\Bigr) \prod_{e'\in F\setminus D} b_{e'}
\prod_{e'\in D\setminus F} 1/b_{e'} = \prod_{e\in F} b_e\ .
$$
The boundary of $F$ coincides with the boundary of $F'$ by the construction of the partial duality. Thus the corresponding arrow monomials in $c$ and $K$'s are also equal to each other. One may check this with the table on page \pageref{tab:cont-del-pd}. 
\end{proof}

\section{Virtual links} \label{s:vir}

Virtual links, introduced in  \cite{Ka3} (a different approach was suggested in  \cite{GPV}), are represented by diagrams similar to ordinary knot diagrams, except some crossings are designated as {\it virtual}. Here are some examples of virtual knots.
$$\risS{-18}{ex}{}{65}{20}{20}\hspace{3cm}
  \risS{-18}{v31}{}{40}{0}{20}\hspace{3cm}
  \risS{-18}{v41}{}{55}{0}{20}
$$

Virtual link diagrams are considered up to plane isotopy, the {\it classical}
Reidemeister moves:
$$\risS{-10}{RI}{}{65}{20}{12}\qquad\qquad
  \risS{-10}{RII}{}{65}{0}{0}\qquad\qquad
  \risS{-10}{RIII}{}{65}{0}{0}\quad ,
$$
and the {\it virtual} Reidemeister moves:
$$\risS{-10}{RI-v}{}{65}{17}{15}\qquad\quad
  \risS{-10}{RII-v}{}{65}{0}{0}\qquad\quad
  \risS{-10}{RIII-v}{}{65}{0}{0}\qquad\quad
  \risS{-10}{RIV-v}{}{65}{0}{0}\quad .
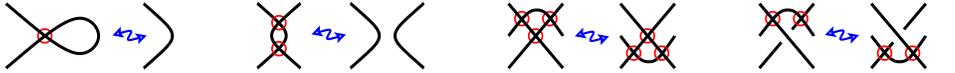
$$

\subsection{Kauffman bracket and Jones polynomial}
The Kauffman bracket for virtual links is defined in the same way as for classical links. Let $L$ be a virtual link diagram.
Consider two ways of resolving a classical crossing.
The {\it $A$-splitting},\ 
$\risS{-4}{cr}{}{15}{15}{8}\ \leadsto\ \risS{-4}{Asp}{}{15}{0}{0},$
is obtained by joining the two vertical angles swept out by the overcrossing arc when
it is rotated counterclockwise toward the undercrossing arc.
Similarly, the {\it $B$-splitting},\ 
$\risS{-4}{cr}{}{15}{15}{8}\ \leadsto\ \risS{-4}{Bsp}{}{15}{0}{0},$
is obtained by joining the other two vertical angles. A {\it state} $s$ of
a link diagram $L$
is a choice of either an $A$ or $B$-splitting at each classical crossing.
Denote by $\cS(L)$ the set of states of $L$.
A diagram $L$ with $n$ crossings has $|\cS(L)| = 2^n$
different states.

Denote by $\a(s)$ and $\b(s)$ the numbers of $A$-splittings and $B$-splittings
in a state $s$, respectively, and by $\d(s)$ the number of
components of the curve obtained from the link
diagram $L$ by 
splitting according to the state $s \in \cS(L)$. Note that virtual crossings do not connect components.

\begin{defn}\label{def:kb}\rm \cite{Ka1}
The {\em Kauffman bracket} of a diagram $L$ is a polynomial in three variables
$A$, $B$, $d$ defined by the formula
$$ \kb{L} (A,B,d)\ :=\ \sum_{s \in \cS(L)} \,
A^{\a(s)} \, B^{\b(s)} \, d^{\,\d(s)-1}\,.
$$
The {\em Jones polynomial}
$J_L(t)$ is obtained from the Kauffman bracket by a simple substitution:
$$A=t^{-1/4},\qquad B=t^{1/4},\qquad d=-t^{1/2}-t^{-1/2}\ ;$$
$$J_L(t)\ := (-1)^{w(L)} t^{3w(L)/4} \kb{L} (t^{-1/4}, t^{1/4}, -t^{1/2}-t^{-1/2})\ ,
$$
where $w(L)$ is the {\it writhe} of the diagram $L$, which is the sum of signs assigned to oriented classical crossings according to the rule:
$$\risS{-10}{or-cr-p}{\put(0,6){$+1$}}{30}{10}{20}\ , \qquad
\risS{-10}{or-cr-m}{\put(0,6){$-1$}}{30}{0}{0}\ .
$$
\end{defn}

Note that $\kb{L}$ is {\em not} a topological invariant of the link; it depends on the link diagram and changes with Reidemeister moves. 

\subsection{Dye-Kauffman arrow polynomial  \cite{DK}}\label{ss:DKarrow}
We can keep more information splitting a classical crossing. Namely, when a splitting does not respect the orientation, we put two arrows on the branches of the splitting oriented counterclockwise near the crossing:
$$\risS{-18}{or-cr-p}{}{40}{20}{20}
\quad\risS{-2}{toto}{}{25}{0}{0}\quad\risS{-18}{or-cr-pa}{}{40}{0}{0}\ ,
\hspace{2cm}\risS{-18}{or-cr-m}{}{40}{20}{20}
\quad\risS{-2}{toto}{}{25}{0}{0}\quad\risS{-18}{or-cr-pa}{}{40}{0}{0}\ .
$$
Thus the state circles are supplied with an arrow structure. With each such circle $c$ we associate the variable $K_c$ as in Definition
\ref{def:aBR}. Then we can define the {\it arrow bracket polynomial} as   
$$ \kb{L}_A (A,B,d)\ :=\ \sum_{s \in \cS(L)} \,
A^{\a(s)} \, B^{\b(s)} \, d^{\,\d(s)-1} \prod_{c\in s} K_{c}\ .
$$
The standard substitution $B:=A^{-1}$, $d:=-A^2-A^{-2}$ gives the {\it normalized Dye-Kauffman arrow polynomial}  \cite{DK}:
$$\langle L\rangle_{NA}:=(-A^3)^{-w(L)}\kb{L}_A (A,A^{-1},-A^2-A^{-2})\ ,
$$
which is an invariant of virtual links. The invariance under the Reidemeister moves follows from the rule of cancellation of arrows in Definition
\ref{def:aBR}. A remarkable observation of H.~Dye and L.~Kauffman is that for classical link diagrams, all arrows will cancel, and the $K$ variables thus do not occur in the arrow polynomial. In this case it is essentially equivalent to the Jones polynomial (after the further substitution $A=t^{-1/4}$).

\section{Arrow Thistlethwaite theorem}

\subsection{From virtual link diagrams to arrow ribbon graphs}

With each state $s$ of a virtual link diagram $L$ we associate an arrow ribbon graph $G_L^s$. The vertices of $G_L^s$ are obtained by gluing discs to the state circles of $s$. The edges of $G_L^s$ correspond to the classical crossings of $L$. Each is obtained by gluing a small planar band connecting the two opposite arcs of the particular splitting of $s$. We consider two types of edge-ribbons.
If a crossing of $L$ is resolved as an $A$-splitting in the state $s$, we assign $+1$ to the corresponding edge, if it is resolved as 
a $B$-splitting, then we assign $-1$. We thus get a signed ribbon graph.

The arrow structure assigns two arrows on opposite sides of each edge-ribbon according to the orientation of the plane where the corresponding planar band is located. 
If a crossing splitting in $s$ respects the orientation of strands, we
put two arrows on the free edge arcs of the corresponding small planar band induced by the counterclockwise orientation of the plane. If the splitting in $s$ of a crossing does not respect the orientation, we put the two arrows on the attaching arcs, again according to the counterclockwise orientation of the plane.

The next example illustrates this construction. 
$$\begin{array}{l}
\risS{-20}{constr1}{\put(-3,40){\mbox{$L$}}\pn{15}{-8}{Diagram}
  \put(113,28){\smiinc}\put(145.5,32){\smiinc}\put(165,15){\smiinc}
      \pn{125}{-8}{State $s$}
  \put(220,28){\smiinc}\put(253.5,33){\smiinc}\put(283,21){\smiinc}
      \pn{210}{-8}{Attaching planar bands}
  \put(337,28){\smiinc}\put(368,34){\smiinc}\put(386,9){\smiinc}
  \pn{335}{-8}{Putting arrows}
  }{400}{40}{90} \\
\hspace{1cm}\risS{-20}{constr2}{
  \pn{24}{-8}{Pulling state circles apart}
  \pn{150}{-8}{Untwisting state circles}
  \put(320,65){\smiinc}\put(355,56){\smiinc}\put(355,8){\smiinc}
  \pn{275}{-8}{Forming the ribbon graph $G_L^s$}
  }{380}{0}{45} 
\end{array}$$

Observe that if we choose the state $\bar{s}$ all of whose splittings do not respect the orientation of the strands, the resulting arrow structure on the ribbon graph $G_L^{\bar{s}}$ coincides exactly with a particular arrow presentation of $G_L^{\bar{s}}$, since all of the arrows are placed on attaching arcs.  

As illustrated in  \cite{Ch}, the partial duals of $G_L^{s}$, for any choice of state $s$, are in bijective correspondence with the states of $L$. A review of the definition of partial duality and our method for assigning arrows to splittings should convince the reader that this correspondence extends to arrow ribbon graphs. We summarize these results in the following lemma. 

\begin{lemma}\label{le:two-st}
Let $s$ and $s'$ be two states of the same diagram $L$. Then the arrow ribbon graphs $G_L^s$ and $G_L^{s'}$ are partial dual with respect to a set of edges corresponding to the crossings where the states $s$ and $s'$ are different from each other. 
\end{lemma}

\begin{thm}[{\bf Arrow Thistlethwaite theorem}]\label{th:av-th}
Let $L$ be a virtual link diagram and let $G_L^s$ be the signed arrow ribbon graph corresponding to a state $s$ with $e_-$ negative edges and $e_+$ positive edges. Then the arrow bracket polynomial of $L$
is a specialization of the arrow dichromatic polynomial of $G_L^s$:
\begin{equation}\label{eq:ath}
\kb{L}_A (A,B,d) = \frac{A^{e_+}B^{e_-}}{d} A_{G_L^s}(1,\bb,d,\bK)
\ ,
\end{equation}
where the weight variables are specialized to \qquad
$b_e=\left\{\begin{array}{ll}
             B/A&\mbox{if $e$ is positive,}\\
             A/B&\mbox{if $e$ is negative.}
           \end{array}\right.
$
\end{thm}

\begin{proof}
Lemma \ref{le:two-st} and Proposition \ref{prop:pardu} imply that the right hand side of \eqref{eq:ath} does not depend on the initial state $s$.

The 1-to-1 correspondence between the states and spanning subgraph of $G_L^s$ is obvious: a state $s'$ corresponds to a spanning subgraph $F$ that contains only the edges corresponding to the crossings of $L$ where $s'$ differs from $s$. It remains to compare the monomials of \eqref{eq:ath} corresponding to $s'$ and $F$. The boundary components of $F$ are the state circles of $s'$, so they carry the same arrow structure, and therefore the $K$ variables are identical in any given monomial. In particular, $bc(F)=\d(s')$, so the exponents of $d$ are also the same.

The exponent of $A$ on the right hand side of \eqref{eq:ath} is equal to 
$$e_+ - e_+(F) + \ e_-(F) = e_+(E(G_L^s)\setminus F)+e_-(F)\ .
$$
where $e_+(F)$ (resp. $e_-(F)$) is the number of positive (resp. negative) edges of the subgraph $F$.
It is easy to see that the last number is equal to $\a(s')$. 

Finally the exponent of $B$ at the right hand side of \eqref{eq:ath} is equal to\qquad 
$e_- + e_+(F) - \ e_-(F) = e_-(E(G_L^s)\setminus F)+e_+(F) = \b(s')\ .
$\end{proof}

\section{Specializations}

It was indicated in  \cite{Ch} that the previously known ribbon graph generalizations of the Thistlethwaite theorem from  \cite{CP,ChVo,DFKLS}
can be unified using different states in the construction of the ribbon graph
$G_L^s$. Here we formulate the corresponding arrow polynomial generalizations.

\subsection{All-$A$-splitting state}
If $s=s_A$ is a state consisting of all $A$-splittings, then all the edges of
$G_L^s$ are positive. In this case all weight variables will be equal to each other: $b_e=B/A$. Theorem \ref{th:av-th} becomes
$$
\kb{L}_A (A,B,d) =  
\sum_{F\subseteq E(G_L^s)}A^{e(\overline{F})} B^{e(F)} d^{bc(F)-1}
   \prod_{f\in\partial(F)}K_{i(f)}\ ,
$$
where $\overline{F}:= E(G_L^s)\setminus F$ is the complementary set of edges.
This equation directly extends the results of  \cite{DFKLS} to the arrow polynomial.

\subsection{Seifert state}
Let $s$ be the Seifert state where all splittings preserve the orientation of the link $L$.
Using \eqref{eq:mBR} we can define an {\it arrow version} of the Bollob\'as-Riordan polynomial as 
$$
ABR_G(X,Y,Z,\bK):=
\sum_{F\subseteq E(G)} \ 
  (\prod_{e\in F}x_e)\ (\prod_{e\not\in F}y_e) 
  X^{r(G)-r(F)}Y^{n(F)}Z^{k(F)-bc(F)+n(F)}
   \prod_{f\in\partial(F)} K_{i(f)}\ .
$$
Substituting\ 
$\displaystyle
x_+=y_+=1,\ x_-=\sqrt{X/Y},\ y_-=\sqrt{Y/X}
$\ 
as in Remark \ref{rem:mBR}, we get the {\it signed unweighted version} of 
the arrow Bollob\'as-Riordan polynomial  
$$
sBR_G(X,Y,Z,\bK)=\sum_{F\subseteq E(G)} \ 
  X^{r(G)-r(F)+s(F)}Y^{n(F)-s(F)}Z^{k(F)-bc(F)+n(F)}
   \prod_{f\in\partial(F)} K_{i(f)}\ ,
$$
where $s(F) := \frac{e_{-}(F)-e_{-}(\wo F)}{2}$.
In this case Theorem \ref{th:av-th} becomes
$$
\kb{L}_A (A,B,d) =  A^{n(G_L^s)} B^{r(G_L^s)} d^{k(G_L^s)-1}
 sBR_{G_L^s}(Ad/B,Bd/A,1/d,\bK)\ ,
$$
which directly extends the results of  \cite{ChVo} to the arrow polynomial.

\end{document}